\documentclass{amsart}
\usepackage{amsmath, amsthm}

    \theoremstyle{plain}
    \newtheorem{Thm}{Theorem}[section]
    \newtheorem{Prop}[Thm]{Proposition}
    \newtheorem{Lemma}[Thm]{Lemma}
    \newtheorem*{Lemma*}{Lemma}
    \newtheorem{Cor}[Thm]{Corollary}

    \theoremstyle{definition}
    \newtheorem{Def}[Thm]{Definition}
    \newtheorem*{Def*}{Definition}
    \newtheorem{Example}[Thm]{Example}

    \theoremstyle{remark}
    \newtheorem{Question}[Thm]{Question}
    \newtheorem{Remark}[Thm]{Remark}
    \newtheorem*{Remark*}{Remark}

    \numberwithin{equation}{section}

    \newcommand{\define}[1]{\emph{#1}}      

    \newcommand{\field}[1]{\mathbb{#1}}
    
    \newcommand{\Z}{\field{Z}}
    
    \newcommand{\R}{\field{R}}
    \newcommand{\C}{\field{C}}

    \DeclareMathSymbol{\fieldk}{\mathalpha}{AMSb}{"7C} 

    \newcommand{\abs}[1]{\lvert#1\rvert}
    
    \newcommand{\norm}[1]{\lVert#1\rVert}

	\newcommand{\trivgp}{\{e\}}
	\newcommand{\genby}[1]{\left<{#1}\right>}
	\newcommand{\half}{{\frac{1}{2}}}

    \DeclareMathOperator{\Tr}{Tr}

    \DeclareMathOperator{\rank}{rank}
    \DeclareMathOperator{\Vol}{Vol}
    \DeclareMathOperator{\Inj}{Inj}
    \DeclareMathOperator{\Diam}{Diam}

    \newcommand{\g}{\gamma}
    \newcommand{\G}{\Gamma}
    \newcommand{\ep}{\varepsilon}
    \newcommand{\la}{\lambda}
    \newcommand{\Laplace}{\Delta}

    \newcommand{\cover}[1]{\widetilde{#1}}
    \newcommand{\torus}{\mathbb{T}}

    \newcommand{\ltwo}{l^{2}}
    \newcommand{\Ltwo}{L^{2}}
    
    \newcommand{\Luck}{L\"{u}ck}

    \DeclareMathOperator{\SH}{short}
    \newcommand{\short}[1]{{\SH(#1)}}
	\newcommand{\Cheb}{Chebyshev}
	\newcommand{\Xp}{X^{\prime}}
    \newcommand{\Gp}{\Gamma^{\prime}}
    \newcommand{\Lp}{\Delta^{\prime}}

\begin{document}

%
%
\title{Growth of Betti Numbers}

\author{Bryan Clair \and Kevin Whyte}
\address{Department of Mathematics,
         Saint Louis University,
         St. Louis, MO 63103}
\email{bryan@slu.edu}
\address{Department of Mathematics,
         University of Illinois at Chicago,
         Chicago, IL 60607}
\email{kwhyte@uic.edu}
\date{November 9, 2001}

\maketitle

%
%
%
\section*{Introduction}

Let $X = \cover{X}/\G$ be a finite simplicial complex.  We study the
growth rate of the Betti numbers of finite regular covers of $X$.  Let
$X_i = \cover{X}/\G_i$, $i = 1,2,\dots$, be a sequence of
finite regular coverings of $X$.
It is easy to see that the sequence of Betti numbers $\{b_q(X_i)\}$
can grow at most linearly in $[\G:\G_i]$.

A theorem of \Luck\ describes exactly when the linear growth rate is achieved,
settling a conjecture of Kahzdan.
\Luck's Theorem~\cite{luck:resfin} states that
\[
   b^{(2)}_{q}(\cover{X};\G) =  \lim_{i \to \infty}
   	                    \frac{b_{q}(X_{i})}{[\G:\G_{i}]}
\]
when the $\G_{i}$ are a tower of finite index normal subgroups of
$\G$ with $\cap \G_i=\trivgp$.

This shows that linear growth of Betti numbers occurs if and only if the
$\Ltwo$-Betti number, $b^{(2)}_q(\cover{X};\G)$, is non zero.
Or, equivalently, if and only if there are non-trivial $\Ltwo$ harmonic
$q$-cochains on $\cover{X}$.
We are interested in the growth rate when the $\Ltwo$-Betti number
vanishes.  Our work follows \Luck's in outline, with more involved
estimates.

There are two key factors in estimating the Betti numbers of covers of
$X$.  The first measures how $X$ is unrolled in the cover:

\begin{Def*}
Let $\G$ be a finitely generated group, with the
word metric from some fixed generating set.  For any subgroup
$\G'$ of $\G$, define $\short{\G'}$ to be the length of the
shortest non-identity element of $\G'$.
\end{Def*}

It is clear from the proof of \Luck's Theorem that instead of a
tower of $\G_{i}$, one could allow any sequence of
$\G_{i} \triangleleft \G$
with $\short{\G_{i}} \to \infty$.

The second ingredient in our estimates is
the behavior of the $\Ltwo$ spectrum of
$\cover{X}$.   Just as \Luck's theorem relates linear growth of Betti
numbers to the kernel of the $\Ltwo$ Laplacian, our
results relate slower growth rates to the behavior of the spectrum of
the $\Ltwo$ Laplacian near zero.

To say ``growth rates'' suggests a tower, or at least a sequence of
coverings.  But the main theorems in this paper are stated as bounds
on given finite covers $\Xp$ of $X$:

\begin{Thm}\label{maintheorem}
Let $X$ be a finite simplicial complex, and $\cover{X}$ an
infinite regular covering with covering group $\G$.
Suppose that $b^{(2)}_{q}(\cover{X};\G) = 0$, or equivalently that
there are no $\Ltwo$ harmonic $q$-cochains on $\cover{X}$.

\begin{enumerate}
\item(Spectral Gap)
Suppose there is a gap near $0$ in the $\Ltwo$
spectrum of $\cover{X}$ in dimension $q$.
Then there are $C > 0$ and $M > 0$ so that for any finite regular
cover $X'=\cover{X}/\Gp$ of $X$:
\[
  b_{q}(\Xp) \leq C\frac{[\G:\Gp]}{ e^{M\short{\Gp}} }.
\]

\item(Positive Novikov-Shubin Invariant) If $\cover{X}$ has Novikov-Shubin
invariant $\alpha_q > 0$, then for any $\ep > 0$ there is a
$C_\ep > 0$ so that for any finite regular cover $X'=\cover{X}/\G'$ of
$X$:
\[
  b_{q}(\Xp) \leq C_{\ep} \frac{[\G:\Gp]}
                         {\bigl(\short{\Gp}\bigr)^{\alpha_{q}-\ep}}.
\]

\item(General Case) For any $X$ there is a $C$ so that for  any finite
regular cover $X'=\cover{X}/\G'$ of $X$:
\[
     b_{q}(\Xp) \leq C \frac{[\G:\Gp]}{\log(\short{\Gp})}.
\]
\end{enumerate}
\end{Thm}

The three cases of this theorem are proved later as
Theorem~\ref{thm:gapgen}, Theorem~\ref{thm:nspos}, and
Theorem~\ref{thm:sublog}, respectively.
They are all special cases of a more general statement
(Proposition~\ref{bettibnd}) which involves the spectrum of $\Laplace$
near $0$ in a more technical way.

For a family of covers $X_{i} = \cover{X}/\G_{i}$ with
$\short{\G_{i}} \to \infty$,
Theorem~\ref{maintheorem} can be interpreted
as giving sublinear upper bounds
on the rate of growth of $b_{q}(X_{i})$ in terms of $[\G:\G_{i}]$.
We can relate $\short{\G_{i}}$ directly to $[\G:\G_{i}]$ if we assume
the covers of $X$ ``unroll'' evenly in all directions, via the
following definition:

\begin{Def*}
A family $\{\G_i\}$ of finite index normal subgroups of $\G$
is \define{uniform} if there is a $C > 0$ so that
\[
  [\G:\G_{i}] \leq \Vol(B_\G(C \cdot \short{\G_{i}}))
\]
for all $i$.
Here $B_{\G}(r)$ is the ball in $\G$ of radius $r$, in the word
metric.

As an example, congruence subgroups of arithmetic groups are uniform.
\end{Def*}

The uniform assumption is most interesting in the case of spectral
gap, where it leads to particularly clean bounds of a form
previously studied by
Sarnak and Xue (\cite{sarnakxue},\cite{mr.x}).

\begin{Thm}\label{gapuniform}
Let $X$ be a finite simplicial complex, and $\cover{X}$ an infinite
regular covering with covering group $\G$.
Suppose that $\G$ has exponential growth, $\{\G_i\}$ is a uniform
family of finite index normal subgroups, and that $\Laplace_{q}$
on $\cover{X}$ has no spectrum below $\la_{0}$.
\begin{itemize}
\item
There is some $C > 0$ and $\beta < 1$ so that
for all $X_{i} = \cover{X}/\G_{i}$,
\[
  b_{q}(X_{i}) < C [\G:\G_{i}]^{\beta}.
\]
\item
Fix $\la < \la_{0}$.
There is some $C > 0$ and $\beta < 1$ so that for all $i$,
\[
   \#\bigl\{ \mu \leq \la \ \big|\
       \text{$\mu$ is an eigenvalue of $\Laplace_{q}$ on $X_{i}$}
       \bigr\}
       <  C [\G:\G_{i}]^{\beta}.
\]
\end{itemize}
\end{Thm}

The second part of Theorem~\ref{gapuniform} is proved by a slight
alteration of the main argument for Betti numbers.  Although one could
make this sort of eigenvalue bound in a more general setting, we
only carry out the calculations in the most interesting case.

To produce geometric applications, one needs a good understanding
of the spectrum of the Laplacian near zero.
The locally symmetric spaces provide a wide class of interesting
examples, as the $\Ltwo$-Betti numbers and Novikov-Shubin invariants are
known.  For these spaces $\short{\Gamma}$ is essentially the injectivity
radius of the cover, and $[\G:\Gp]$ is the volume of the cover.
Specifically, for hyperbolic manifolds:

\begin{Thm}\label{thm:hyper}
Let $X$ be a compact hyperbolic $n$ manifold.  There are constants
$C > 0$ and $\beta_{q} > 0$
so that for any finite regular cover $X'$ of $X$ one has
\begin{itemize}
 \item For $n$ odd:
   \begin{itemize}
   \item If $q \neq \frac{n\pm 1}{2}$, then
         \[ b_q(X') \leq C \frac{\Vol(X')}{e^{\beta_q \Inj(X')}}. \]
   \item For $q=\frac{n\pm 1}{2}$,
         \[ b_q(X') \leq C \frac{\Vol(X')\cdot\log{\Inj(X')}}{\Inj(X')}. \]
   \end{itemize}
 \item For $n$ even:
   \begin{itemize}
   \item If $q \neq \frac{n}{2}$, then
         \[ b_q(X') \leq C \frac{\Vol(X')}{e^{\beta_q \Inj(X')}}. \]
   \item For $q=\frac{n}{2}$, the $\Ltwo$-Betti number
         $b_{q}^{(2)}(\mathbb{H}^{n}; \G)$ is non-zero.
   \end{itemize}
\end{itemize}
\end{Thm}

We investigate the case $\G=\Z^n$ in depth, using different
techniques:

\begin{Thm}
Suppose $\cover{X}$ is a regular $\Z^{n}$ covering of a finite
simplicial complex $X$, and assume the $\Ltwo$-Betti number
$b^{(2)}_{q}(\cover{X};\Z^{n}) = 0$.
Then there is a constant $C > 0$ so that
for any finite cover $X'=\cover{X}/\G'$ of $X$:
\[
  b_{q}(\Xp) \leq C \frac{[\Z^{n}:\Gp]}{\short{\Gp}}.
\]
\end{Thm}

This is similar to the bound Theorem~\ref{maintheorem} would give if
$\alpha_q(X)$, the Novikov-Shubin invariants of $X$, satisfied $\alpha_q
\geq 1$.  However, Lott \cite[Ex. 42]{lott:heat} gives examples
with $\G = \Z$ and $\alpha_{q}(X)$ arbitrarily small.

%
%
\section{Stripes}
\label{stripes}

To understand the growth of Betti numbers of covers it helps to have
some examples.  In particular, one would like a method for constructing
complexes with a given fundamental group where it is easy to describe
finite covers and their Betti numbers.  The construction
described in this section starts with an arbitrary finite complex
and attaches pieces which do not change
the fundamental group and which contribute in a straightforward way
to the Betti numbers of finite covers.

Given complexes $X$ and $Y$, and a complex $Z$ which includes as a
subcomplex of both $X$ and $Y$, build a complex $X \cup_{Z} Y$
which is the quotient of the disjoint union of $X$ and $Y$ by the
relation which identifies the copies of $Z$.  More generally, given
a $Z$ with maps $f:Z \to X$ and $g: Z \to Y$, build a new complex,
which we describe as ``gluing $X$ and $Y$ along $Z$'', by carrying out the
previous construction where we replace $X$ and $Y$ by the mapping
cylinders of $f$ and $g$, which are homotopy equivalent to $X$ and $Y$
and contain $Z$ as a subcomplex.

Assume that $\pi_{1}(X)=\G$.  By Van Kampen's theorem,
$\pi_{1}(X \cup_{Z} Y)$ is also equal to $\G$ provided that
$\pi_{1}(Z) \hookrightarrow \pi_{1}(X)$ and $\pi_{1}(Z)=\pi_{1}(Y)$.
With these assumptions on fundamental groups, put
\[
   W = X \cup_{Z} Y.
\]

Let $\Gp$ be a finite index normal subgroup of $\G$,
and let $W'$ be the cover of $W$ with $\pi_{1}(W') = \Gp$.
We can describe
$W'$ as the result of a gluing $X'\cup_{Z'}Y'$ where $X'$, $Y'$, and
$Z'$ are the induced finite covers of $X$, $Y$, and $Z$.

To understand the homology of $W'$, we need to
understand the homology of $X'$, $Y'$, and $Z'$, as well as
the maps on homology induced by the inclusions.  The special case
we want to consider here is where all of the homology comes from
the parts we have glued on.  So assume that $X$ has dimension $n$,
and that we want to compute $H_{q}(W')$ for $q>n$.  This homology is
simply the homology of $Y'$.

If $\pi_{1}(Z)=\G$ then we have
gained nothing by this construction, as it reduces to the same
problem for one of the pieces.  However, if $\pi_{1}(Z)$ is a proper
subgroup of $\G$ then the examples coming from this construction
give useful insight into the general problem.

\begin{Example}
As a simple example, suppose $\pi_{1}(Z) = \trivgp$.
Then $Y'$ is the disjoint union of $[\G:\Gp]$ copies of $Y$.
If $H_{q}(Y)=0$, then $H_{q}(W') = 0$ as well, for any $\Gp$.
If $H_{q}(Y)\neq 0$, then the Betti numbers $b_{q}(W')$
grow linearly in $[\G:\Gp]$.
Neither case is interesting from the point of view of our problem.
\end{Example}

In general, the number of components of $Y'$ is
\[
  b_{0}(Y') = [\pi_{1}(Y) : \pi_{1}(Y)\cap\Gp],
\]
each of which is the cover of
$Y$ with fundamental group $\pi_{1}(Y) \cap \Gp$.
Thus, to
understand the Betti numbers of these covers, we need to understand
the the Betti numbers of the covers of $Y$ and the number of
components.  The key source of examples is the following:
\begin{Def}
Put $Z=S^{1} \hookrightarrow X$, and $Y=S^{1} \times S^{q}$.
The construction of $W = X \cup_{Z} Y$ is called
\define{gluing a stripe to $X$ along $Z$}.
We will often refer to $Z$ simply by the corresponding
element of $\G = \pi_{1}(X)$.
\end{Def}

For the stripe construction, the components of the cover $Y'$ of $Y$
are $S^{1} \times S^{q}$, so their $q$-dimensional homology has rank
one.  Thus
\[
  \dim(H_{q}(W')) = [\pi_{1}(Z) : \pi_{1}(Z)\cap\Gp].
\]
Let $\g$ be a generator of $\pi_{1}(Z) = \Z$, and let $o(\g)$ denote
the order of $\g$ in the quotient $\G/\Gp$.  Then
\[
  \dim(H_{q}(W')) = \frac{[\G:\Gp]}{o(\g)}.
\]

To get examples of fast sublinear growth of Betti numbers, we want
$o(\g)$ to grow, but as slowly as possible.   The following Lemma
demonstrates the key piece of geometry of $\Gp$ that enters our
bounds:

\begin{Lemma}\label{lem:stripe}
Suppose $\g \in \G$ has infinite order.  There is
a $C>0$ so that for any complex $W$ built by gluing a stripe
along $\g$, and for any finite index normal $\Gp$:
\[
  \dim(H_{q}(W')) \leq C\frac{[\G:\Gp]}{\short{\Gp}}.
\]
\end{Lemma}

\begin{proof}
By the earlier discussion, this amounts to the claim that the order
of $\g$ in $\G/\Gp$ is at least
$\short{\Gp}/C$.   If $\g^{k} = e$ in $\G/\Gp$
then $\g^{k} \in \Gp$.  Since $\g$ has infinite order
$\g^{k}$ is a non-trivial element of $\Gp$, and so has
length at least $\short{\Gp}$.  On the other hand, the length of
$\g^{k}$ is at most $k$ times the length of $\g$.  Thus the
Lemma follows with $C$ equal to the length of $\g$.
\end{proof}


\begin{Example}\label{ex:tight}
Let $\G = \Z^{2} = \genby{a}\times\genby{b}$, and take a family of
subgroups
$\G_{i} = m_{i}\Z\times n_{i}\Z$.  Put $X = \torus^{2}$, and glue
a $q$-stripe to $X$ along $a$.  Then $b_{q}(W_{i}) = n_{i}$, while
$[\G:\G_{i}] = m_{i}n_{i}$.

If $n_{i} > m_{i}$, then $\short{\G_{i}} = m_{i}$ and the bound in
Lemma~\ref{lem:stripe} is achieved.  If $m_{i} > n_{i}$, however, it
is clear that the stripe along $a$ fails to capture the geometry.
A stripe along $b$ would give a larger Betti number.
\end{Example}

In the previous example, we could have glued two stripes so that
for any fixed family of groups $\G_{i} = m_{i}\Z\times n_{i}\Z$,
growth of order $[\G:\G_{i}]/\short{\G_{i}}$ is achieved.
However, for some families $\G_{i}$, there is no stripe
construction that achieves this rate of growth:
\begin{Example}\label{ex:nogrow}
Let $\G=\Z^{2}$,
let $\{p_{k}\}$ be the sequence of primes,  and let $\G_{i}$ be
the kernel of the map to $\Z/p_{1}p_{3}\dotsm p_{2i+1} \times
\Z/p_{2}p_{4}\dotsm p_{2i}$.  Then, $\short{\G_{i}}=p_{2}\dotsm
p_{2i}$.
For any $\g = (n,m) \in \Z^{2}$, the index of
$\genby{\g}$ in $\G/\G_{i}$ is bounded by $nm$.  Therefore, the Betti
numbers of the $\G_{i}$ covers of the complex obtained by gluing
a stripe (or any finite number of stripes) are bounded over all $i$.
\end{Example}

We will prove that this example is not a defect of the stripe
construction.  Indeed, for $\G=\Z^{n}$, and any family of
coverings, the stripes examples achieve the fastest possible growth of
Betti numbers (see Remark~\ref{rem:stripe}).
We do not know if this is true for general $\G$.

%
%
\section{The Abelian Case}
\label{abelian}
In this section, we investigate in depth the special case
where $\G = \Z^{n}$.  The techniques are different from the general
situation considered in Section~\ref{growthbnds}. Here, the key
step is an algebraic estimate of rational solutions to a
trigonometric polynomial.  The results are of a similar form
to the general theorems, and illustrate the strength of the stripe
construction.

\begin{Thm}\label{thm:Zn}
Suppose $\cover{X}$ is a regular $\Z^{n}$ covering of a finite
simplicial complex $X$, and assume the $\Ltwo$-Betti number
$b^{(2)}_{q}(\cover{X};\Z^{n}) = 0$.
Then there is a constant $C(X)$ so that
for any $\Gp < \Z^{n}$, we have
\begin{equation}\label{eq:zbound}
  b_{q}(\Xp) \leq C(X)\frac{[\Z^{n}:\Gp]}{\short{\Gp}}.
\end{equation}
Here $\Xp = \cover{X}/\Gp$.
\end{Thm}

Suppose that $X$ has $a$ cells in dimension $q$.
All work is done in the fixed dimension
$q$, but the $q$ is suppressed in most of the notation.
Put $\G = \Z^{n}$.

Choosing lifts of cells from $X$ to $\cover{X}$, we identify
$C^{(2)}_{q}(\cover{X})$ with $\bigoplus_{1}^{a}\ltwo(\G)$.  In this
basis, $\Laplace$ is represented by an $a\times a$ matrix $B$ with
entries in $\Z[\G]$.

Since $\G = \Z^{n}$ is abelian, we can form $\det(B) \in \Z[\G]$.
The determinant $\det(B)$ will not depend on the choice of lifts of cells,
so we will simply write $\det(\Laplace)$.

Identify the $n$-torus $\torus^{n}$ with the irreducible unitary
representations of $\Z^{n}$.  Let
\[
  K = \{ \rho \in \torus^{n}\ |\ \rho(\det(\Laplace)) = 0 \}.
\]
Call $K$ the \define{pattern} of $X$ on $\torus^{n}$.

We want to calculate the Betti numbers for the finite cover $\Xp$ of
$X$.  The Laplacian on $q$-chains on $\Xp$ is also given by the
matrix $B$, now acting on $\bigoplus_{1}^{a}\ltwo(\G/\Gp)$.
The space $\ltwo(\G/\Gp)$ splits as a direct sum of irreducible
representations $\rho$, where $\rho$ lies on a lattice $\Lambda$ of
rational points in $\torus^{n}$.  $\Lambda$ is (non-canonically)
isomorphic to $\G/\Gp$, and
\[
  \Lambda = \left\{ \rho \in \torus^{n}\ |\
                   \rho\ \text{is trivial on}\ \Gp \right\}.
\]

\begin{Lemma}\label{TorusEstimate}
\[
  \abs{\Lambda \cap K} \leq b(\Xp) \leq a \cdot \abs{\Lambda \cap K}
\]
\end{Lemma}
\begin{proof}
We have
\[
  \bigoplus_{1}^{a}\ltwo(\G/\Gp) =
     \bigoplus_{\rho\in\Lambda}
       \left( \bigoplus_{1}^{a}V_{\rho} \right),
\]
where $V_{\rho} = \C$.  The Laplacian $\Lp$ preserves the
$a$-dimensional space $\bigoplus_{1}^{a}V_{\rho}$,
and has kernel there precisely when $\rho(\det(\Laplace)) = 0$.
The dimension of the kernel is obviously bounded by $a$.
\end{proof}

The goal, then, is to count certain rational points in $\torus^{n}$
which lie on the pattern of $X$.  The pattern of $X$ is the set
of root-of-unity solutions to a polynomial equation, and these
are fully described by Conway and Jones~\cite{conwayjones}.

Choose generators $g_{1},\dotsc,g_{n}$ for
$\Z^{n}$.  Then $\det(\Laplace)$ is a polynomial in the $g_{k}$
(with integer coefficients).
For $\rho = (x_{1},\dotsc,x_{n}) \in \torus^{n}$,
$\rho(g_{k}) = e^{2\pi i x_{k}}$.  Therefore,
\[
  \rho(det(\Laplace)) = \sum_{I}A_{I}e^{2\pi i l_{I}(x)},
\]
where $l_{I}(x)$ is some integer linear combination of the $x_{i}$.
It follows from~\cite{conwayjones} that the rational
solutions to
\[
  \sum_{I}A_{I}e^{2\pi i l_{I}(x)} = 0
\]
lie on one of a finite family of rational linear subspaces
of $\torus_{n}$.

We have shown that the rational points in the pattern of $X$
are contained in a finite collection of rational linear subspaces.
Fix one such space $L \subset \torus^{n}$.
If $L$ is not itself a subgroup, then it is a coset of some subgroup
$L^{\prime} < \torus^{n}$, and then
\[
  \abs{\Lambda \cap L} \leq \abs{\Lambda \cap L^{\prime}}.
\]
So, to estimate $\abs{\Lambda \cap L}$ we may as well assume
that $L$ is a subgroup of $\torus^{n}$.

There is a subgroup $\hat{L} \subset \G$ so that
\[
  L = \left\{ \rho \in \torus^{n}\ |\
                   \rho\ \text{is trivial on}\ \hat{L} \right\}.
\]
Note that $\dim(L) + \rank_{\Z}\hat{L} = n$.

Choose any $g \in \hat{L}$.  For optimal constants, one should
make $\norm{g}$ as small as possible.
Let $o(g)$ denote the order of $g$ in $\G/\Gp$.  We have
\[
   o(g) \geq \frac{\short{\Gp}}{\norm{g}}.
\]

Construct a representation
$\rho:\G \to S^{1}$ which is trivial on $\Gp$ and sends $g$ to
a primitive $o(g)$-th root of unity (this is easy, since $\G/\Gp$ is
a finite abelian group).
For $k = 1,2,\dotsc,o(g)$, each representation $\rho^{k}$ is in
$\Lambda$, but all assume different values on $g$.  In particular,
$\rho^{1},\rho^{2},\dotsc,\rho^{o(g)}$ all lie in different cosets
of $\Lambda \cap L < \Lambda$.  Therefore
\[
   \abs{\Lambda \cap L} \leq
      \frac{\abs{\Lambda}}{o(g)} \leq
      \frac{\norm{g}\cdot[\G:\Gp]}{\short{\Gp}}.
\]

Summing over the various $L$ which appear in $K$ (the pattern of $X$),
we conclude that
\begin{equation}\label{LambdaCapK}
   \abs{\Lambda \cap K} \leq
      \frac{C_{1}(X)\cdot[\G:\Gp]}{\short{\Gp}}.
\end{equation}
Here $C_{1}(X)$ is a constant depending on $X$, or more precisely,
depending on the rational linear subspaces appearing in the pattern
of $X$.

Combining \eqref{LambdaCapK} with Lemma~\ref{TorusEstimate} completes
the proof of Theorem~\ref{thm:Zn}.
\qed

As a corollary of the proof, we can make an even better statement
for $\G = \Z$:

\begin{Cor}
Suppose $\cover{X}$ is a regular $\Z$ covering of a finite simplicial
complex $X$.  Put $X_{i} = \cover{X}/i\Z$, the $i$-fold covering of
$X$.  Then exactly one of the following possibilities occurs:
\begin{enumerate}
\item $b^{(2)}_{q}(\cover{X};\Z) \neq 0$, so the sequence of
      Betti numbers $b_{q}(X_{i})$ is asymptotically linear in $i$.
\item There is a constant $C(X)$ so that $b_{q}(X_{i}) \leq C(X)$
      for all $i$.
\end{enumerate}
\end{Cor}
\begin{proof}
The torus $\torus^{1}$ is just $S^{1}$, so a linear subspace is either
a single point or the entire circle.  If the pattern of $X$ contains
the entire circle, we are in the first case.  Otherwise, the pattern
of $X$ consists of a finite collection of $k$ points (counted with
multiplicities), and so the
betti numbers of $X_{i}$ are bounded by $k$ times the number of cells
in $X$.
\end{proof}

\begin{Remark}\label{rem:stripe}
For $\G = \Z^{n}$ and for a given
sequence of subgroups $\G_{i} \subset \G$, we need only consider
the stripe construction when looking for large Betti numbers.
This is because given $X = \cover{X}/\G$, the Betti numbers $b_{q}(X_{i})$
depend only on the rational linear subspaces contained in the
pattern of $X$.  We can reproduce these rational linear subspaces
by repeatedly gluing stripes to an $n$-torus.
\end{Remark}

\begin{Remark}
Example~\ref{ex:tight} demonstrated a sequence $\G_{i} \subset \Z^{n}$
and a space $X$ for which the bound in Theorem~\ref{thm:Zn} is tight.
However, it is not hard to construct a sequence of $\G_{i}$ for which
the bound in Theorem~\ref{thm:Zn} is far from the best possible.
See Example~\ref{ex:nogrow}.
\end{Remark}

\begin{Remark}
The $q^{th}$ Novikov-Shubin invariant for a space with a $q$-stripe
is $1$, the same as for $\R$.  Compare the bound
\eqref{eq:zbound} for $\G$ abelian with the general bound
\eqref{eq:nsbound} for Novikov-Shubin invariant $1$.
\end{Remark}

%
%
\section{General Bounds}
\label{growthbnds}
In this section, we derive a general bound on Betti numbers of
coverings, in terms of the $\Ltwo$ spectral density function.
The argument follows the proof of \Luck's Theorem
\cite{luck:resfin}, but we carefully control the estimates
throughout.

\subsection{Preliminaries}
Let $X$ be a finite simplicial complex, and $\cover{X}$ an infinite
regular cover, with covering transformation group $\G$.
Suppose we have a normal subgroup $\Gp$ of finite index in $\G$.
Form $\Xp = \cover{X}/\Gp$, a covering of $X$ of order $[\G:\Gp]$.

Let $\Laplace$ be the combinatorial Laplacian on
$C^{(2)}_{q}(\cover{X})$, the space of $q$-dimensional
$\ltwo$-cochains.  Let $\Lp$ be the
Laplacian on $q$-cochains of $\Xp$.

Suppose $X$ has $a$ cells in dimension $q$.  Lifting these cells to
$\cover{X}$ gives a basis over $\ltwo(\G)$ for
$C^{(2)}_{q}(\cover{X})$.
These lifts descend to give a basis
over $\ltwo(\G/\Gp)$ for $C_{q}(\Xp)$.

In this basis, $\Laplace$ is represented by
an $a\times a$ matrix $B$ with entries in $\Z[\G]$, acting by right
multiplication on $\bigoplus_{j=1}^{a}\ltwo(\G)$.
The same matrix $B$ describes $\Lp$,
now acting by right multiplication on
$\bigoplus_{j=1}^{a}\ltwo(\G/\Gp)$.

\begin{Lemma}\label{lem:Kbound}
There exists a number $K > 1$ so that $\norm{\Laplace} \leq K$
and so that for any group $\Gp$ as above, $\norm{\Lp} \leq K$.
\end{Lemma}
\begin{proof}
This is identical to \cite[Lemma 2.7]{luck:resfin}.  The number $K$
can be defined from the coefficients of group elements appearing in $B$.
\end{proof}

\begin{Lemma}\label{lem:poly}
Let
\[
  R = \max\bigl\{\norm{g}\ \big|
            \ \text{$g\in\G$ appears with nonzero coefficient in some
            $B_{ij}$} \bigr\}.
\]
Then for any polynomial $p$ with
$\deg(p) < \frac{\short{\Gp}}{R}$,
\[
  \Tr_{\G}p(\Laplace) = \Tr_{\G/\Gp}p(\Lp)
\]
\end{Lemma}
\begin{proof}
Write
  \[ \sum_{j=1}^{a}(p(B))_{jj} =
     \sum_{g\in\G}\lambda_{g}g \quad (\in \Z[\G]). \]
Then
  \[ \Tr_{\G}p(\Laplace) = \lambda_{e} \]
and
  \[ \Tr_{\G/\Gp}p(\Lp) = \sum_{g\in\Gp}\lambda_{g}. \]
However, $B$ contains group elements of length at most $R$, so $B^{n}$
contains group elements of length at most $nR$.  Therefore $p(B)$
contains group elements of length less than $\short{\Gp}$, and so
$\lambda_{g} = 0$ for all $g \in \Gp$, $g \neq e$.
\end{proof}

\subsection{Density functions and Betti numbers}
\newcommand{\Pn}{{\mathcal{P}_{n}}}
\begin{Def}
Let
\[
  \Pn = \left\{ \quad\text{Polynomials $p$ in one variable}\quad
        \left|
        \begin{array}{c}
           \deg p \leq n;\\
           \text{$p$ non-negative on $[0,1]$};\\
           p(0) = 1
        \end{array}
        \right.\right\}
\]

For a probability measure $d\mu$ on $[0,1]$, put
\[
  J(n,\mu) = \inf_{p \in \Pn}\int_{0}^{1}p(x)d\mu(x).
\]
\end{Def}

The measure of interest to us is the spectral density function
of the Laplacian, suitably rescaled.  To define it,
let $\{E(\la) : \la \in [0,\infty)\}$ denote the
family of spectral projections of $\Laplace$.  Then
the spectral density function of $\Laplace$ is
\[
   F(\la) = \Tr_{\G}E(\la).
\]

\begin{Prop}\label{bettibnd}
Let
\[ n < \frac{\short{\Gp}}{R} \]
and
\[ \mu(x) = \frac{F(Kx)}{a}. \]
Then
\[
  b_{q}(\Xp) \leq a[\G:\Gp]J(n,\mu).
\]
\end{Prop}
\begin{proof}
Let $p \in \Pn$.  For any $r > 1$, there is some $\la > 0$ so that
every $x \in [0,\la]$ satisfies
\[
  rp(\frac{x}{K}) \geq 1.
\]
Let $\chi$ be the characteristic function of the interval $[0,\la]$.
We have $rp(\frac{x}{K}) \geq \chi(x)$ on $[0,K]$.
Then
\begin{align*}
  \frac{b_{q}(\Xp)}{[\G:\Gp]} &\leq
     \Tr_{\G/\Gp}\chi(\Laplace)\\
     &\leq \Tr_{\G/\Gp} rp(\frac{\Laplace}{K})\\
     &= r \Tr_{\G} p(\frac{\Laplace}{K})
\end{align*}
using Lemma~\ref{lem:poly}.  Sending $r \to 1$,
\begin{align*}
  b_{q}(\Xp) &\leq
    [\G:\Gp]\Tr_{\G}p(\frac{\Laplace}{K})\\
    &= [\G:\Gp]\int_{0}^{K}p(\frac{\la}{K}) dF(\la)\\
    &= a[\G:\Gp]\int_{0}^{1}p(x)d\mu(x).
\end{align*}
Since $p \in \Pn$ was arbitrary, the proof is done.
\end{proof}

\subsection{Choosing a polynomial}\label{sec:poly}
We now proceed to estimate $J(n,\mu)$, with an eye towards
using the decay of $\mu(x)$ as $x \to 0$.
Fix $z \in (0,1)$.  Let $p \in \Pn$ be bounded by 1 on
$[0,z]$ and bounded by $M$ on $[z,1]$.   If $p(z) < M$, the
estimates will be strictly improved by moving $z$ to the left, so we
may as well restrict attention to the situation where $p(z)$ is the
maximal value of $p$ on $[z,1]$.  We have
\begin{align*}
  J(n,\mu) &\leq \int_{0}^{1}p(x)d\mu(x) \\
           &= \int_{0}^{z}p(x)d\mu(x) + \int_{z}^{1}p(x)d\mu(x) \\
           &\leq \mu(z) + p(z).
\end{align*}

To minimize $\mu(z) + p(z)$, we want a polynomial which drops as
quickly as possible from $x=0$ to $x=z$ and then stays low until $x=1$.
Via a linear transformation, this is equivalent to finding a
polynomial which is bounded by $\pm 1$ on $[-1,1]$ and
grows as quickly as possible for $x > 1$.  Among polynomials of
degree $n$, the \Cheb\ polynomial $T_{n}$ is the optimal solution
\cite{rivlin}.

The \Cheb\ polynomials are defined by $T_{n}(\cos(\theta)) = \cos(n\theta)$.
The first few are $1, x, 2x^{2}-1, 4x^{3}-3x$.  We need two facts:  first,
$T_{n}(1) = 1$ for all $n$;  second, the \Cheb\ polynomials satisfy a
recurrence relation
\[ T_{n+1}(x) = 2xT_{n}(x) - T_{n-1}(x) \]
which yields
\[
  T_{n}(x) = \half\left[\left(x + \sqrt{x^{2}-1}\right)^{n} +
                \left(x - \sqrt{x^{2}-1}\right)^{n}\right].
\]

\newcommand{\ratiopm}[1]{{\frac{1+#1}{1-#1}}}
Now, put
\[
  l(x) = \left(\frac{-2}{1-z}\right)x + \left(\ratiopm{z}\right),
\]
so $l(1) = -1$, $l(z) = 1$, and $l(0) = \ratiopm{z}$.  Set
\[
  p_{n}(x) = \frac{T_{n}\bigl(l(x)\bigr) + 1}{T_{n}\bigl(l(0)\bigr) + 1}.
\]
It is easy to see that $p_{n} \in \Pn$.  Therefore
\begin{align*}
  J(n,\mu) &\leq \mu(z) + p_{n}(z) \\
           &\leq \mu(z) + \frac{2}{T_{n}\left(\ratiopm{z}\right)}\\
           &= \mu(z) + \frac{4}
                            {{\left(\ratiopm{\sqrt{z}}\right)}^{n} +
                             {\left(\ratiopm{\sqrt{z}}\right)}^{-n}}\\
           &\leq \mu(z) + 4{\left(\ratiopm{\sqrt{z}}\right)}^{-n}\\
           &\leq \mu(z) + 4e^{-2n\sqrt{z}}.
\end{align*}
This proves:
\begin{Prop}
For any $n > 0$, probability measure $d\mu$, and $z \in [0,1]$,
\begin{equation}\label{Jbound}
  J(n, \mu) \leq \mu(z) + 4e^{-2n\sqrt{z}}
\end{equation}
\end{Prop}

To choose $z$ appropriately requires information about $\mu$.
We handle three specific cases in the next sections.

%
%
\section{Spectral Gap}
\newcommand{\lao}{\lambda_{0}}
Suppose that the $\Ltwo$ spectrum for $\cover{X}$
has a gap at zero.  This means that there exists
$\lao > 0$ with $F(\la) = 0$ for all $\la < \lao$.

Continuing notation from Section~\ref{growthbnds},
$\mu(x) = F(Kx)/a$ has a gap at zero of size $\lao/K$.
Letting $z \to \lao/K$ in \eqref{Jbound} yields
\[
  J(n,\mu) \leq 4e^{-2n\sqrt{\lao/K}}.
\]
Combining with Proposition~\ref{bettibnd} proves the following:
\begin{Thm}\label{thm:gapgen}
If $\cover{X}$ has spectral gap of size $\lao$
in dimension $q$, then
\begin{equation}\label{gapshort}
  b_{q}(\Xp) \leq 4a\frac{[\G:\Gp]}{ e^{M\short{\Gp}} },
\end{equation}
where $M = \frac{2}{R}\sqrt{\frac{\lao}{K}}$.
Here $R$ is the maximum
length of group elements appearing in the matrix for $\Laplace$,
and $K$ is a bound on the norms of $\Laplace$ and
$\Laplace_{i}$ as in Lemma~\ref{lem:Kbound}.
\end{Thm}


If $\G$ has exponential growth, $[\G:\Gp]$ will be
at least exponential in $\short{\Gp}$.
For a family of subgroups $\G_{i} \triangleleft \G$,
a natural assumption is that $\short{\G_{i}}$
actually does grow like $\log[\G:\G_{i}]$.  This is the assumption
of the following Corollary.

\begin{Cor}
Let $\{\G_{i}\}$ be a family of normal, finite index subgroups
of $\G$.  Suppose there exists $D > 0$ so that
\begin{equation}\label{loggrowth}
  \short{\G_{i}} > D\log[\G:\G_{i}] - \text{constant}.
\end{equation}
Put $X_{i} = \cover{X}/\G_{i}$.
Then there is some $C > 0$ so that
\[
  b_{q}(X_{i}) < C \cdot [\G:\G_{i}]^{1-MD}
\]
for all $i$.
\end{Cor}

This also proves the first part of Theorem~\ref{gapuniform}.  The
estimate \eqref{loggrowth} comes directly from the
assumption that $\G$ has exponential growth and that the family
$\{\G_{i}\}$ is uniform.
Applying \eqref{loggrowth} to the next Proposition proves the second
part of Theorem~\ref{gapuniform}.

\begin{Prop}
Fix $\la < \lao$, and let
\[
  E(\la) =    \bigl\{ \mu \leq \la \ \big|\
       \text{$\mu$ is an eigenvalue of $\Lp_{q}$ on $X'$}
       \bigr\}.
\]
Then there is a $C > 0$ and $M > 0$ so that
\[
   \#E(\la) < C\frac{[\G:\Gp]}{ e^{M\short{\Gp}} }.
\]
\end{Prop}
\begin{proof}
We continue the notation of Section~\ref{sec:poly}, where $p_{n}$ is a
linear transformation of the $n^{th}$ \Cheb\ polynomial.  Note that
$p_{n}$ is decreasing on $[0,\lao/K]$.  Then
\[
   \#E(\la) \cdot p_{n}(\la/K)
      \leq \sum_{\mu \in E(\la)} p_{n}(\mu/K)
      \leq [\G:\Gp] \Tr_{\G/\Gp} p_{n}(\Lp/K).
\]
If $n < \short{\Gp}/R$, Lemma~\ref{lem:poly} applies, and we get
\[
   \#E(\la) \leq [\G:\Gp] \frac{\Tr_{\G}p_{n}(\Laplace/K)}{p_{n}(\la/K)}.
\]
The spectrum of $\Laplace$ is empty below $\lao$, so
\[
   \#E(\la) \leq [\G:\Gp] \frac{p_{n}(\lao/K)}{p_{n}(\la/K)}.
\]
Finally, a simple calculation shows that the ratio
$\frac{p_{n}(\lao/K)}{p_{n}(\la/K)}$ decays exponentially in $n$ as
$n \to \infty$ (although the base approaches 1 when $\la$ nears $\lao$).
Replacing $n$ with $\short{\Gp}/R$ completes the proof.
\end{proof}

%
%
\section{Positive Novikov-Shubin Invariant}
\label{NSsection}
The $q^{th}$ Novikov-Shubin invariant of $\cover{X}$ describes the
decay of the spectral density function $F(\la)$ near $\la = 0$.
Suppose there is some $C > 0$ so that
\begin{equation}\label{eq:ns}
   C^{-1}\lambda^{\alpha_{q}/2} < F(\la) < C\lambda^{\alpha_{q}/2}
\end{equation}
for small $\la$.  Then $\alpha_{q}$ is the $q^{th}$
Novikov-Shubin invariant of $\cover{X}$.
In general, one defines
\[
  \alpha_{q} = 2 \liminf_{\la\to 0^{+}}
               \frac{\log(F(\la) - F(0))}{\log(\la)}
             \quad \in [0,\infty]
\]
when $F$ has no spectral gap.  The $2$ is a normalization which is
discarded by some authors.  With our definition, all Novikov-Shubin
invariants of $\R^{n}$ are $n$.

For the moment, suppose that $\mu(x) < Cx^{\beta}$
for some $C > 0$ and $\beta > 0$.

To make $J(n,\mu)$ small, we choose $z$ so that the two terms in the
bound \eqref{Jbound} are roughly the same size for large $n$.
We ask that
\[
  z^{\beta} = e^{-2n\sqrt{z}}
\]
or
\[
  \beta\log(z) = -2n\sqrt{z}.
\]
Substituting $x = z^{-1/2}$ yields
\[
  x\log(x) = n/\beta,
\]
which has the solution $\log(x) = W(n/\beta)$, where $W$ is the
Lambert W-function.  For large $n$, $W(y)$ is asymptotic to
$\log(y/\log(y))$.  Based on this, we choose
\[
  z = {\left(\frac{\log(n/\beta)}{n/\beta}\right)}^{2}.
\]

Now plug in to the bound \eqref{Jbound} for $J(n,\mu)$:
\begin{align*}
  J(n,\mu) &\leq C{\left(\frac{\log(n/\beta)}{n/\beta}
                        \right)}^{2\beta} +
                4\exp\left(-2n\frac{\log(n/\beta)}{n/\beta}
                     \right)\\
           &= C{\left(\frac{\log(n/\beta)}{n/\beta}
                        \right)}^{2\beta} +
                4\left(\frac{1}{n/\beta}\right)^{2\beta}.
\end{align*}
Therefore, there is some constant $C^{\prime}$ depending on $C$ and
$\beta$ so that for all $n$,
\[
  J(n,\mu) \leq C^{\prime} \cdot
               {\left(\frac{\log(n)}{n}\right)}^{2\beta}.
\]

Combining this result with Proposition~\ref{bettibnd} gives the
following:
\begin{Thm}
Suppose the $q^{th}$ spectral density function of $\cover{X}$ satisfies
\[ F_{q}(\la) < C\la^{\beta} \]
for some $\beta > 0$, $C > 0$.
Then there is a constant $C_{1} > 0$ so that
\begin{equation}\label{eq:nsbound1}
  b_{q}(\Xp) \leq C_{1} [\G:\Gp]
      {\left(\frac{\log(\short{\Gp})}{\short{\Gp}}\right)}
      ^{2\beta}.
\end{equation}
\end{Thm}

Now, we interpret this for Novikov-Shubin invariants:
\begin{Thm}\label{thm:nspos}
Suppose $b^{(2)}_{q}(\cover{X};\G) = 0$ and
the Novikov-Shubin invariant $\alpha_{q}(\cover{X})$ is
positive.  Then for every $\ep > 0$, there is some constant $C_{\ep}$
depending on $\ep$ and $X$ so that for any finite covering
$\Xp = \cover{X}/\Gp$ of $X$,
\begin{equation}\label{eq:nsbound}
  b_{q}(\Xp) \leq C_{\ep} \frac{[\G:\Gp]}
                         {\bigl(\short{\Gp}\bigr)^{\alpha_{q}(\cover{X})-\ep}}.
\end{equation}
\end{Thm}
\begin{proof}
If $\alpha_{q}(\cover{X}) > 0$, we have a bound of the form
\[
  F_{q}(\la) < C\la^{\beta}
\]
for any $\beta < \alpha_{q}(\cover{X})/2$.
Since we can always make $\beta$ a little larger than needed for
a given $\ep$, the $\log(\short{\Gp})$ term in \eqref{eq:nsbound1}
is absorbed into the constant.
\end{proof}

%
%
\section{Sublogarithmic Decay}
The most general estimate known for spectral density functions is
\[
   F(\la) < \frac{a\log(K)}{-\log(\la)}.
\]
This estimate holds whenever $\G$ belongs a large class of
groups $\mathcal{G}$.
The class $\mathcal{G}$ includes all residually finite groups.
For details, see \cite{dm:amen}, \cite{clair:resamen}, \cite{schick:approx}.

Assume that $\mu(x) < \frac{C}{-\log(x)}$.
As in Section~\ref{NSsection},
we want to choose $z$ so that the two terms in the
bound \eqref{Jbound} are roughly the same size for large $n$.
We ask that
\[
  \frac{1}{-\log(z)} = e^{-2n\sqrt{z}}.
\]
Substitute $x = 2n\sqrt{z}$ to get
\[
  e^{x} + 2\log(x) = 2\log(2n),
\]
which has the approximate solution $x = \log(2\log(2n))$ for large $n$.
As the precise constants won't matter in the end, choose
\[
  z = {\left(\frac{\log(\log n)}{n}\right)}^{2}.
\]

Now plug in to the bound \eqref{Jbound} for $J(n,\mu)$:
\[
  J(n,\mu) \leq \frac{C}{-2\bigl(\log(\log(\log n)) - \log(n)\bigr)}
                 + 4e^{-\log(\log n)}
\]
Therefore, there is some constant $C^{\prime}$ so that for all $n$,
\begin{equation}\label{sublogbnd}
  J(n,\mu) \leq \frac{C^{\prime}}{\log(n)}.
\end{equation}

\begin{Thm}\label{thm:sublog}
Given $\G$ and $X = \cover{X}/\G$, with $b^{(2)}_{q}(\cover{X};\G) = 0$.
There is a constant $C$ so that for any finite
covering $\Xp = \cover{X}/\Gp$ of $X$,
\[
     b_{q}(\Xp) \leq C \frac{[\G:\Gp]}{\log(\short{\Gp})}.
\]
\end{Thm}
\begin{proof}
If $\G$ is not residually finite, $\short{\Gp}$ is uniformly
bounded over all $\Gp$, so the theorem is vacuously true.
Otherwise, from \cite{luck:resfin}, one has decay of the $\Ltwo$
spectrum of $\cover{X}$ of the form
\[
   F(\la) < \frac{a\log(K)}{-\log(\la)},
\]
for $\la < \ep$, with some $\ep > 0$.
Now $\mu(x) < \frac{C}{-\log(x)}$ for some $C > 0$, and we apply
\eqref{sublogbnd} to Proposition~\ref{bettibnd}.
\end{proof}

\begin{Remark}
If $\G$ is not residually finite, form $\G_{f} = \bigcap \Gp$, where
the intersection is over all finite index normal subgroups of $G$.
Then all finite covers of $X$ are covered by $X_{f} = \cover{X}/\G_{f}$.
The group $\G/\G_{f}$ acts on $X_{f}$ and is residually finite, so we
can still bound Betti numbers of $\Xp$ in a non-trivial way.
The subtle point is that $\short{\Gp}$ is replaced by the shortest
element of $\Gp/\G_{f}$ as a subgroup of $\G/\G_{f}$.
\end{Remark}

\section{Applications And Questions}

In applying Theorem~\ref{maintheorem} to specific spaces there are two
main obstacles: controlling the $\Ltwo$ spectral data and calculating the
growth of $\short{\G'}$ for various covers.

First, if $X$ is a non-positively curved complex and $\G=\pi_1(X)$ then
$\short{\G'}$ is comparable to the injectivity radius $\Inj(X')$.
In particular, one has the upper bound
\[
  \Inj(X') < \Diam(X').
\]
One can often estimate $\Diam(X')$ in terms of volume,
which is multiplicative for covers.
(If $\G$ is a quotient of $\pi_1(X)$, then
$\short{\G'}$ is comparable to the ``injectivity radius of the projection
$\cover{X} \to X'$'', meaning the supremum of $r>0$ so that $r$-balls in
$\cover{X}$ are embedded in $X$.  This is equivalent to half the length of
the shortest closed geodesic in $X'$ which does not lift to a loop in
$\cover{X}$.)

This estimate for $\short{\G'}$ has a parallel for arbitrary
groups:
\begin{equation}\label{eqshortdiam}
   \short{\G'} \leq 1+\Diam(\G/\G')
\end{equation}
If one knows the growth rate
of the volume of balls in $\G$, \eqref{eqshortdiam}
gives an upper bound for
$\short{\G'}$ in terms of $[\G:\G']$.  In particular, if $\G$ has
exponential growth then one knows that for some $C$ and any
$\G'$, $\short{\G'} \leq C \log{[\G:\G']}$.

Unfortunately, we are usually interested in lower bounds for
$\short{\G'}$, not upper bounds.  If the covers $X'$ of $X$ unroll
in all directions equally, then these upper bounds are accurate
estimates.

\begin{Def}
A family $\{\G_i\}$ of finite index normal subgroups of $\G$
is \define{uniform} if there is a $C > 0$ so that
\[
  [\G:\G_{i}] \leq \Vol(B_\G(C \cdot \short{\G_{i}}))
\]
for all $i$.
Here $B_{\G}(r)$ is the ball in $\G$ of radius $r$, in the word
metric.

Say that a family of regular covers $\{X_{i}\}$ is uniform
if the corresponding family of groups is uniform.
\end{Def}

\begin{Lemma}
If $\G$ is arithmetic, the collection of congruence
subgroups $\{\G_{n}\}$ is uniform.
\end{Lemma}
\begin{proof}
We show this when $\G$ has exponential growth.  The proof for
polynomial growth is no more difficult.  View $\G$ as a group of integral matrices.
Let $m$ be the largest matrix entry appearing in all generators of $\G$.
A word in generators of $\G$ must be of length at least
$\log_{m}(n)$ to have an entry of size $n$.  Therefore,
$\short{\G_{n}} \geq \log_{m}(n)$.
As the volume growth of balls in $\G$ is exponential, we can
choose a $C$ so that $\Vol(B_\G(C \cdot \short{\G_{n}}))$
is larger than any given polynomial in $n$, and in particular is larger
than the index $[\G:\G_{n}]$.
\end{proof}

\begin{Example}
Suppose $X$ is an $n$-dimensional manifold of pinched negative curvature,
and $\cover{X}$ is the universal cover of $X$.
$\cover{X}$ has spectral gap for $\abs{q - n/2} \geq 1$ \cite{dx:spec}.
For such a $q$, one has $C > 0$ and $M > 0$ so that for all
finite regular covers $X'$ of $X$:
\[
    b_{q}(X') < C \frac{\Vol(X')}{e^{M \Inj(X')}}
\]
For a uniform family of covers $\{X_{i}\}$, Theorem~\ref{gapuniform}
applies and gives $C > 0$, $\beta<1$ so that
\begin{equation}\label{eq:volpow}
   b_q(X_{i}) < C \Vol(X_{i})^\beta
\end{equation}
for all $i$.
\end{Example}

In particular, \eqref{eq:volpow} holds for congruence covers of arithmetic
hyperbolic manifolds, outside the middle dimension.
Here, Sarnak and Xue~\cite{sarnakxue} have conjectured an upper bound of
the form \eqref{eq:volpow}, with $\beta = 2q/(n-1) + \ep$.
For these manifolds, lower bounds of the same type are known.
Xue~\cite{mr.x} has shown:
\[
    b_{q}(X_{p}) > C_{\ep} \cdot \Vol(X_{p})^{\delta-\ep}
\]
for any $\ep > 0$, and with an explicit $\delta$.

For more general symmetric spaces, the spectrum near zero is well
understood.  The $\Ltwo$ Betti numbers vanish except possibly in the
middle dimension.
Near the middle dimension, the spectral density obeys a power law
decay, and away from the middle range there is spectral gap.
See~\cite{olbrich} for the exact results.
This allows us to give upper bounds for the growth
rates of Betti numbers of locally symmetric manifolds, for example
Theorem \ref{thm:hyper}.

Adams and Sarnak~\cite{adamssarnak}
have precise results for similar questions, calculating
multiplicities of representations instead of Betti numbers.

Notice that in
\begin{Question} Is there a $C>0$ so that for any compact hyperbolic
3-manifold $M$,
  \[  b_1(M) \leq C \frac{\Vol(M)}{\Inj(M)}?  \]
More generally, for a compact locally symmetric space
modeled on $G/K$,
how big can the Betti numbers be in terms of volume and injectivity radius?
\end{Question}

For more general spaces, computing the $\Ltwo$ spectral data is very hard.
Even the question of whether all spaces have positive Novikov-Shubin
invariants is open.  If true, our bounds for arbitrary spaces would get
significantly better, changing from
\[
   \frac{[\G:\G']}{\log{\short{\G'}}}
\] to \[
   \frac{[\G:\G']}{{\short{\G'}}^\beta}
\]
(for some $\beta$ depending on $X$).
Indeed, our results could be used to detect a space with vanishing
Novikov-Shubin invariant:
\begin{Question}
Is there a space $X$ and sequence of finite normal covers $X_{i} =
\cover{X}/\G_{i}$ so that the $q^{th}$
Betti numbers of $X_{i}$ grow sublinearly in
$[\G:\G_{i}]$, but faster than
\[
   \frac{[\G:\G_{i}]}{\bigl(\short{\G_{i}}\bigr)^{\alpha}}.
\]
for all $\alpha > 0$?
\end{Question}

We know of no example with Betti numbers growing anywhere near that
quickly:

\begin{Question}
Is there a space $X$ and sequence of finite normal covers $X_{i} =
\cover{X}/\G_{i}$ so that the $q^{th}$ Betti numbers of
$X_{i}$ grow sublinearly in
$[\G:\G_{i}]$, but faster than
\[
   \frac{[\G:\G_{i}]}{\short{\G_i}}.
\]
\end{Question}

Such a space would have $q^{th}$ Novikov-Shubin invariant at most one.
Conceivably, there could even be such examples among hyperbolic
$3$-manifolds.

It would be interesting to translate these questions and the results
of this paper into the setting of extended $\Ltwo$ cohomology,
introduced by Farber~\cite{far:extl2}.

%
%
\bibliographystyle{plain}

\end{document}